%% file: main.tex
\pdfoutput=1
\documentclass{article}
\usepackage{amsmath, amsthm, amssymb}
\usepackage{eucal}
\usepackage{mathrsfs}
\usepackage{geometry}
\usepackage[bottom]{footmisc}
\usepackage[unicode,bookmarksnumbered=true,colorlinks=true,allcolors=blue,linktoc=all]{hyperref}
\usepackage[capitalise,nameinlink]{cleveref}
\usepackage{tikz-cd}
\usepackage{stmaryrd}
\usepackage{listings}
\usepackage{etoolbox}
\usepackage{thm-restate}
\usepackage[nottoc,notlot,notlof]{tocbibind}
\usepackage{enumitem}
\usepackage{mathtools}
\usepackage{quiver}
\usepackage{caption}

\crefname{subsection}{Subsection}{subsections}

\hypersetup{bookmarksdepth=2}
\newtheorem{theorem}{Theorem}

\newtheorem{thm}{Theorem}[section]
\newtheorem{lem}[thm]{Lemma}
\newtheorem{prop}[thm]{Proposition}

\newtheorem{question}[thm]{Question}

\theoremstyle{definition}
\newtheorem{defn}[thm]{Definition}

\newtheorem{assum}[thm]{Assumption}

\newtheorem{remark}[thm]{Remark}

\newcommand{\ncmd}{\newcommand}

\definecolor{DefColor}{rgb}{0.6,0.15,0.25}
\newcommand{\mdef}[1]{\textcolor{DefColor}{#1}}
\newcommand{\tdef}[1]{\mdef{\emph{#1}}}

\ncmd{\mbb}[1]{\mathbb{#1}}
\ncmd{\mrm}[1]{\mathrm{#1}}
\ncmd{\mcl}[1]{\mathcal{#1}}
\ncmd{\mfk}[1]{\mathfrak{#1}}
\ncmd{\mbf}[1]{\mathbf{#1}}
\ncmd{\mscr}[1]{\mathscr{#1}}

\ncmd{\todo}[1]{\textbf{TODO #1}}
\ncmd{\reftodo}[1]{\textbf{REF #1}}

\DeclareRobustCommand{\minwidthbox}[2]{%
  \mathmakebox[\ifdim#2<\width\width\else#2\fi]{#1}%
}
\ncmd{\too}[1][]{\xrightarrow{\minwidthbox{#1}{1em}}}
\ncmd{\iso}{\too[\smash{\raisebox{-0.5ex}{\ensuremath{\scriptstyle\sim}}}]}
\ncmd{\oot}[1][]{\xleftarrow{\minwidthbox{#1}{1em}}}
\ncmd{\osi}{\oot[\smash{\raisebox{-0.5ex}{\ensuremath{\scriptstyle\sim}}}]}
\ncmd{\hooktoo}[1][]{\xhookrightarrow{\minwidthbox{#1}{1em}}}
\ncmd{\adj}[1][]{\mathrel{\substack{\xrightarrow{\minwidthbox{#1}{1em}} \\[-.7ex] \xleftarrow{\minwidthbox{#1}{1em}}}}}

\ncmd{\qin}{\quad\in\quad}

\ncmd{\Id}{\mrm{Id}}
\ncmd{\Nm}{\mathrm{Nm}}

\ncmd{\BB}{\mrm{B}}
\ncmd{\clB}{\mcl{B}}
\ncmd{\CC}{\mcl{C}}
\ncmd{\DD}{\mcl{D}}
\ncmd{\EE}{\mcl{E}}
\ncmd{\TT}{\mcl{T}}
\ncmd{\MM}{\mcl{M}}
\ncmd{\KK}{\mrm{K}}
\ncmd{\KU}{\mrm{KU}}
\ncmd{\cX}{\mcl{X}}
\ncmd{\cY}{\mcl{Y}}
\ncmd{\Kn}{\KK(n)}
\ncmd{\Knm}{\KK(n{-}1)}
\ncmd{\Knp}{\KK(n{+}1)}
\ncmd{\Ko}{\KK(1)}
\ncmd{\Tn}{\mrm{T}(n)}
\ncmd{\Tnp}{\mrm{T}(n{+}1)}
\ncmd{\To}{\mrm{T}(1)}
\ncmd{\KTnp}{\KK_{\Tnp}}
\ncmd{\KKo}{\KK_{\Ko}}
\ncmd{\KTo}{\KK_{\To}}
\ncmd{\THH}{\mrm{THH}}
\ncmd{\TC}{\mrm{TC}}
\ncmd{\TCm}{\TC^-}
\ncmd{\rmE}{\mrm{E}}
\ncmd{\En}{\rmE_n}
\ncmd{\Enp}{\rmE_{n{+}1}}
\ncmd{\JWnp}{\rmE(n{+}1)}
\ncmd{\cJWnp}{\widehat{\JWnp}}
\ncmd{\Enpk}{\Enp(\kappa)}
\ncmd{\lsF}{\mscr{F}}
\ncmd{\lsG}{\mscr{G}}
\ncmd{\bbC}{\mbb{C}}
\ncmd{\GG}{\mbb{G}}
\ncmd{\NN}{\mbb{N}}
\ncmd{\ZZ}{\mbb{Z}}
\ncmd{\QQ}{\mbb{Q}}
\ncmd{\clQ}{\mcl{Q}}
\ncmd{\Qab}{\QQ(\zeta_\infty)}
\ncmd{\Zp}{\ZZ_p}
\ncmd{\Qp}{\QQ_p}
\ncmd{\Fp}{\mbb{F}_p}
\ncmd{\Fpbar}{\overline{\mbb{F}}_p}
\renewcommand{\SS}{\mbb{S}}
\ncmd{\SKn}{\SS_{\Kn}}
\ncmd{\SKnp}{\SS_{\Knp}}
\ncmd{\SKo}{\SS_{\Ko}}
\ncmd{\STn}{\SS_{\Tn}}
\ncmd{\STnp}{\SS_{\Tnp}}
\ncmd{\rmP}{\mrm{P}}
\ncmd{\OO}{\mcl{O}}
\ncmd{\one}{\mbf{1}}

\ncmd{\irchi}[2]{\raisebox{\depth/2}{$#1\chi$}}
\DeclareRobustCommand{\rchi}{{\mathpalette\irchi\relax}}
\ncmd{\ch}{\rchi}
\ncmd{\cch}{\widehat{\scalebox{1.15}{$\rchi$}}}

\ncmd{\fL}{L}
\ncmd{\fLLam}{\fL^\Lambda}
\ncmd{\LKn}{L_{\Kn}}
\ncmd{\LKnp}{L_{\Knp}}
\ncmd{\LKo}{L_{\Ko}}
\ncmd{\LTnp}{L_{\Tnp}}
\ncmd{\Ln}{L_n}
\ncmd{\Lnm}{L_{n{-}1}}
\ncmd{\Mn}{M_n}
\ncmd{\Mnf}{M_n^f}
\ncmd{\Lnf}{L_n^f}
\ncmd{\Lnmf}{L_{n{-}1}^f}
\ncmd{\Mz}{M_0}
\ncmd{\Lz}{L_0}

\ncmd{\LL}{\mrm{L}}
\ncmd{\RR}{\mrm{R}}
\ncmd{\BC}{\mrm{BC}}
\ncmd{\psa}{\oplus}
\ncmd{\dbl}{\mrm{dbl}}
\ncmd{\op}{\mrm{op}}
\ncmd{\pifin}{\pi\text{-}\mathrm{fin}}
\ncmd{\pfin}{p\text{-}\mathrm{fin}}
\ncmd{\seg}{\mrm{seg}}
\ncmd{\perf}{\mrm{perf}}

\ncmd{\pt}{\mrm{pt}}
\ncmd{\Perf}{\mrm{Perf}}
\ncmd{\Mod}{\mrm{Mod}}
\ncmd{\cMod}{\widehat{\Mod}\vphantom{\Mod}}
\ncmd{\Vect}{\mrm{Vect}}
\ncmd{\Ab}{\mrm{Ab}}
\ncmd{\Fin}{\mrm{Fin}}
\ncmd{\Spaces}{\mcl{S}}
\ncmd{\Spacespifin}{\Spaces_{\pifin}}
\ncmd{\Spacespfin}{\Spaces_{\pfin}}
\ncmd{\Sp}{\mrm{Sp}}
\ncmd{\SpTn}{\Sp_{\Tn}}
\ncmd{\SpTnp}{\Sp_{\Tnp}}
\ncmd{\SpKn}{\Sp_{\Kn}}
\ncmd{\SpKo}{\Sp_{\Ko}}
\ncmd{\SpKnp}{\Sp_{\Knp}}
\ncmd{\Span}{\mrm{Span}}
\ncmd{\Spano}{\Span_1}
\ncmd{\Cat}{\mrm{Cat}}
\ncmd{\Catpifin}{\Cat_{\pifin}}
\ncmd{\Catpfin}{\Cat_{\pfin}}
\ncmd{\CatLn}{\Cat_{\Ln}}
\ncmd{\CatLnm}{\Cat_{\Lnm}}
\ncmd{\CatMn}{\Cat_{\Mn}}
\ncmd{\CatMnf}{\Cat_{\Mnf}}
\ncmd{\CatLnf}{\Cat_{\Lnf}}
\ncmd{\Catperf}{\Cat_{\perf}}
\ncmd{\PrL}{\mrm{Pr}^\mrm{L}}
\ncmd{\PrLst}{\PrL_{\mrm{st}}}
\ncmd{\PrLstw}{\PrL_{\mrm{st},\omega}}
\ncmd{\PrLKn}{\PrL_{\Kn}}
\ncmd{\PrLKnw}{\PrL_{\Kn,\omega}}
\ncmd{\Mfg}{\mcl{M}_\mrm{fg}}

\ncmd{\yon}{\text{\usefont{U}{min}{m}{n}\symbol{'110}}}
\DeclareFontFamily{U}{min}{}
\DeclareFontShape{U}{min}{m}{n}{<-> dmjhira}{}

\DeclareMathOperator{\Ind}{Ind}

\DeclareMathOperator{\CAlg}{CAlg}
\DeclareMathOperator{\CMon}{CMon}
\DeclareMathOperator{\coCMon}{coCMon}
\ncmd{\CMoninf}{\CMon_\infty}
\ncmd{\coCMoninf}{\coCMon_\infty}
\DeclareMathOperator*{\colim}{colim}

\ncmd\noloc{%
  \nobreak
  \mspace{6mu plus 1mu}
  {:}
  \nonscript\mkern-\thinmuskip
  \mathpunct{}
  \mspace{2mu}
}

\title{Chromatic higher semiadditivity via redshift}
\author{Shay Ben-Moshe\thanks{Max Planck Institute for Mathematics, Bonn, Germany.} \thanks{Faculty of Mathematics and Computer Science, Weizmann Institute of Science, Israel.}}
\date{}

\begin{document}
	\maketitle
	
	\begin{abstract}
		We give a new proof of the $\infty$-semiadditivity of $\Kn$-local spectra.
		The proof proceeds by induction on the height via algebraic K-theory, utilizing recent advances in the redshift conjecture, instead of using the Ravenel--Wilson computation of the Morava K-theory of Eilenberg--MacLane spaces.
		Along the way, we prove the higher semiadditivity of $\Tn$-local modules over the $\Tn$-localized K-theory of $\Knm$-local ring spectra, which suffices for the applications of higher semiadditivity to the telescope conjecture.
	\end{abstract}
	

	\vspace{.7em}

	\begin{figure}[ht!]
	  \centering
	  \includegraphics[width=140mm]{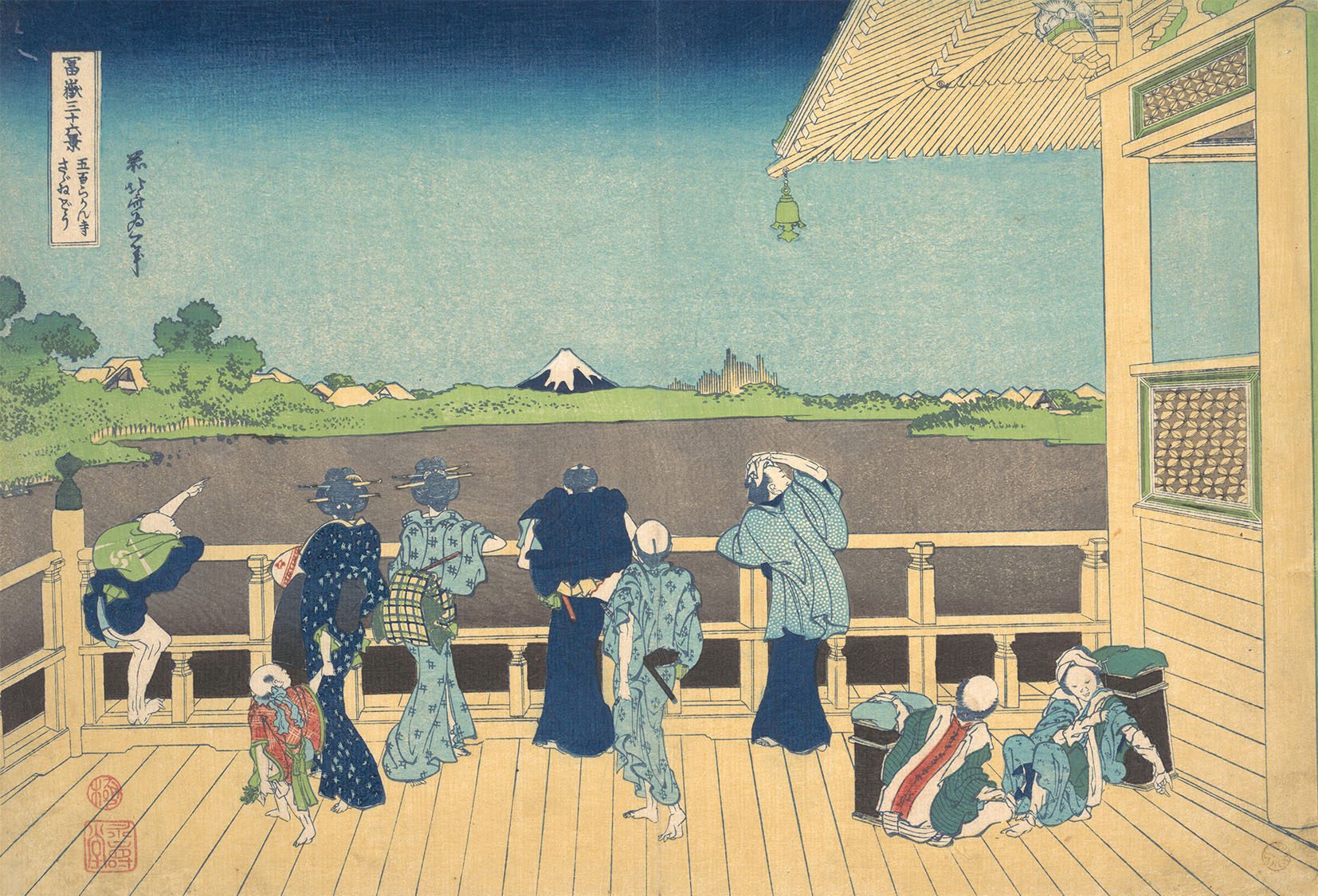}
	  \caption*{
		Sazai Hall at the Temple of the Five Hundred Arhats, by Katsushika Hokusai.
	  }
	\end{figure}

	\pagebreak
	
	\tableofcontents
	
	\input{intro.tex}

\input{height-1.tex}
	\input{proof.tex}

	\bibliographystyle{alpha}
	\bibliography{refs}

\end{document}

%% file: intro.tex
\section{Introduction}

\subsection*{Background}

In their remarkable work on higher semiadditivity, Hopkins--Lurie \cite{HL} showed that the category\footnote{Throughout this paper, we use the term category to mean an $\infty$-category.} of $\Kn$-local spectra is $\infty$-semiadditive, extending the $\Kn$-local Tate vanishing result of Greenlees--Hovey--Sadofsky \cite{GS,HSa}.
Loosely speaking, this says that $\pi$-finite spaces satisfy Poincar\'e duality as though they were $0$-dimensional oriented manifolds, identifying homology and cohomology with no twist or shift.
More precisely, they have shown that for any local system $X\colon A \to \SpKn$ indexed by a $\pi$-finite space $A$, the (inductively defined) norm map
\[
	\Nm\colon \colim_A X \iso \lim_A X
	\qin \SpKn
\]
is an isomorphism.
At the core of their proof is a delicate analysis of the case of the constant $\BB^m C_p$-shaped diagram on the height $n$ Lubin--Tate spectrum $\En$.
In this case, they prove that the norm map gives a perfect pairing on homology by a refinement of the foundational calculations made by Ravenel--Wilson \cite{RW}.
The general case in turn follows by a sequence of formal reductions.

Subsequently, Carmeli--Schlank--Yanovski \cite{TeleAmbi} extended this result to show that $\SpTn$ is $\infty$-semiadditive (which in particular gives a new proof for $\Kn$-local spectra).
Their proof proceeds by induction on the truncatedness of the space $A$, using Kuhn's $\Tn$-local Tate vanishing \cite{Kuhn} as the base of the induction.
To prove $(m{+}1)$-semiadditivity, they show that it suffices to find an $m$-finite $p$-space $A$ with $\pi_m(A) \neq 0$, whose cardinality $|A|_{\STn}$ is invertible, which by the nilpotence theorem is equivalent to the invertibility of $|A|_{\En}$.
To find such a space, they start with $\BB^m C_p$, whose cardinality at $\En$ does not vanish by the Ravenel--Wilson calculation, and use a clever construction to transform it into a suitable space with invertible cardinality.

In our joint work with Carmeli, Schlank and Yanovski \cite{Desc}, we have used this result to prove the higher descent theorem, showing that $\Tnp$-localized K-theory commutes with $\pi$-finite $p$-space limits and colimits of $\Lnf$-local categories, using \cite{DescVan} for $1$-finite spaces as a base for an induction.
In a subsequent paper \cite{card}, we have applied this, together with redshift and the chromatic nullstellensatz due to Yuan and Burklund--Schlank--Yuan \cite{yuan,null}, to compute the cardinality $|A|_{\En}$ for all $\pi$-finite $p$-spaces $A$.
In fact, as we explained in \cite[Observation 10]{card}, this allows one to circumvent the full strength of the Ravenel--Wilson calculation in the proof of the $\infty$-semiadditivity of $\SpTn$.
Indeed, to prove $(m{+}1)$-semiadditivity one only needs to know the non-vanishing of $|\BB^k C_p|_{\En}$ for $k \leq m$, which can be computed in this way assuming only $m$-semiadditivity.

\subsection*{Overview}

In this paper we give a new proof of $\infty$-semiadditivity of $\Kn$-local spectra.
The proof proceeds by induction on the height, passing information to the next height using categorification, or more precisely, using algebraic K-theory and the redshift philosophy, without appealing to the Ravenel--Wilson computation.
Assuming inductively that $\SpKn$ is $\infty$-semiadditive we prove the following.

\begin{theorem}[{\cref{main-thm}}]\label{main-thm-intro}
	The category $\SpKnp$ is $\infty$-semiadditive.
\end{theorem}

The key intermediate step is showing that the category of $\Tnp$-local modules over the $\Tnp$-localized K-theory of $\Kn$-local ring spectra is $\infty$-semiadditive.
We note that this level of generality suffices for the applications of higher semiadditivity to the Burklund--Hahn--Levy--Schlank disproof of the telescope conjecture \cite{tel}.

\begin{theorem}[{\cref{KTnp-sa}}]\label{KTnp-sa-intro}
	The category $\Mod_{\LTnp\KK(\SKn)}(\SpTnp)$ is $\infty$-semiadditive.
\end{theorem}

We now sketch the proof in more detail.
As in \cite{TeleAmbi}, we assume that $\SpKnp$ is already $1$-semiadditive,\footnote{Our argument depends on $1$-semiadditivity implicitly via the descent and purity theorems \cite{DescVan,purity}.} for which a short proof has appeared \cite{Tate}.

The first tool we use is the category $\CatMn$ of $n$-monochromatic categories, namely those stable categories whose mapping spectra are $\Ln$-local and $\Lnm$-acyclic.
This category is interesting for two different reasons.
First, the $\Tnp$-localized K-theory of $\Ln$-local categories factors through $\CatMn$ by the purity theorem \cite{purity,DescVan} as we explained in \cite[Proposition 2.24]{Desc}.
Second, the inductive assumption that $\SpKn$ is $\infty$-semiadditive implies that $\CatMn$ is $\infty$-semiadditive, as we have shown in \cite[Theorem 2.29]{Desc}.

We then use the following key observation regarding the aforementioned higher descent theorem for $\Tnp$-localized K-theory.
The proof in \cite{Desc} only used $\infty$-semiadditivity to reduce from arbitrary limits and colimits along $\pi$-finite $p$-spaces to the case of constant colimits.
Thus, dropping the $\infty$-semiadditivity assumption, we can still prove that $\Tnp$-localized K-theory preserves \emph{constant} $\pi$-finite $p$-space colimits.\footnote{There is a subtlety in the case $n=0$, see \cref{height-1} for an extended discussion.}
We also note that $\Tnp$-localized K-theory is multiplicative, i.e.\ lax symmetric monoidal.

Recall that ($p$-typical) $\infty$-semiadditivity is equivalent to the condition that, for every $\pi$-finite $p$-space, the pairing on its homology (induced by the norm map) is perfect, i.e.\ a self-duality datum.
Since $\Tnp$-localized K-theory is multiplicative and preserves constant $\pi$-finite $p$-space colimits, it carries the relevant self-duality data from $\CatMn$ to $\SpTnp$.
Applying this to the unit, we get the self-duality over $\LTnp\KK(\SKn)$, which proves \cref{KTnp-sa-intro}.

We then use Yuan's proof of redshift \cite{yuan} to transfer this self-duality to the $S^1$-Tate construction on a Morava E-theory of height $n{+}1$.\footnote{Alternatively, one can use the chromatic nullstellensatz \cite{null} to directly pass to a Morava E-theory.}
Finally, we descend to the $\Knp$-local sphere, as in the formal argument of Hopkins--Lurie mentioned above, which concludes the proof of \cref{main-thm-intro}.

\pagebreak

We close the overview with a question.

\begin{question}
	Can the proof be extended to show inductively that $\SpTnp$ is $\infty$-semiadditive?
\end{question}

\begin{remark}\label{telescopic-remark}
	The proof works as-is telescopically until the final descent step.
	More precisely, assuming inductively the $\infty$-semiadditivity of $\SpTn$, the proof shows that we have self-duality over $\LTnp(E^{tS^1})$ for any $E \in \CAlg(\SpTnp)$, for instance over $\LTnp(\STnp^{tS^1})$.
	The remaining step is to descend from $\LTnp(E^{tS^1})$ to the $\Tnp$-local sphere.
	We do not know how to carry out this step, which in the $\Knp$-local case crucially uses our understanding of the $S^1$-Tate construction on complex oriented ring spectra.
\end{remark}

\subsection*{Acknowledgements}

I thank Maxime Ramzi, Tyler Lawson and Asaf Yekutieli for helpful discussions.
I also thank Lior Yanovski, Tomer Schlank and Shachar Carmeli for useful comments on a draft of this paper.
I thank Shachar Carmeli for informing me that Akhil Mathew independently developed a similar argument.
Finally, I would like to thank the Max Planck--Weizmann joint postdoctoral program for its financial support.

%% file: height-1.tex
\section{Height 1}\label{height-1}

There is a subtlety not addressed in the explanation in the introduction in the case of passing from height $n=0$ to height $1$.
The proof of the higher descent theorem relies on the fact that $\Tnp$-localized K-theory preserves sifted colimits, which only holds for $n\geq1$.
To overcome this, in \cite{Desc} we have used the following fact.

\begin{prop}\label{height-1-trunc}
	Let $X$ be a $\Ko$-local spectrum and let $A$ be a simply connected $\pi$-finite $p$-space.
	Then the fold map
	\[
		X[A] \too X
		\qin \SpKo
	\] 
	is an isomorphism.
\end{prop}

Unfortunately, in \emph{op.\ cit.}\ the proof of this fact relies on the $\infty$-semiadditivity of $\SpKo$, which would put us in a circular situation.
Luckily, this result is essentially an immediate consequence of the computation of the complex K-theory of Eilenberg--MacLane spectra due to Anderson--Hodgkin \cite{AH} (of which the Ravenel--Wilson computation is a generalization).

As another alternative, we note that the techniques of \cite{AmbiHeight} give a very short proof of the $\infty$-semiadditivity of $\SpKo$ (assuming $1$-semiadditivity).
Indeed, in this case it is easy to show that $|\BB C_p|$ is invertible, which implies $\infty$-semiadditivity without resorting to the more sophisticated argument from \cite{TeleAmbi} used for higher heights (where one only knows that the cardinalities of Eilenberg--MacLane spectra are non-zero), as we now explain.

\begin{prop}\label{1-inf-sa}
	The category $\SpKo$ is $\infty$-semiadditive.
\end{prop}

\begin{proof}
	Recall from \cite[Corollary 3.3.10]{TeleAmbi} that
	\[
		\dim(\KU_p^{\BB C_p}) = |\BB C_p|_{\KU_p} |C_p|_{\KU_p}
		\qin \pi_0(\KU_p) \simeq \Zp.
	\]
	Since $\dim(\KU_p^{\BB C_p}) = p$ and $|C_p|_{\KU_p} = p$, we conclude that $|\BB C_p|_{\KU_p} = 1$.
	Next, note that the map $\pi_0(\SKo) \to \pi_0(\KU_p)$ exhibits the target as the quotient of the source by the nilradical, and in particular $|\BB C_p|_{\SKo}$ is invertible.
	The result thus follows from \cite[Theorem 3.2.7]{AmbiHeight}.
\end{proof}

\begin{remark}
	The cardinality $|\BB C_p|_{\SKo}$ was explicitly computed in \cite{CY}.
\end{remark}

With this in mind, \cite[Proposition 3.2.3(2)]{AmbiHeight} implies \cref{height-1-trunc}.

%% file: proof.tex
\section{Characterizations of higher semiadditivity}

As was mentioned in the introduction, to prove the $\infty$-semiadditivity of a category, one needs to check that the homology of $\pi$-finite spaces is self-dual.
When the category in question is symmetric monoidal, and the tensor product commutes with $\pi$-finite colimits, it suffices to check that the unit is self-dual in a suitably multiplicative way.
This idea was made precise, and encapsulated in a technically useful way, in the work of Harpaz \cite{Harpaz}.
Let $\Span(\Spacespifin)$ denote the symmetric monoidal category of spans of $\pi$-finite spaces, where the tensor product is given by the cartesian product in the underlying category $\Spacespifin$.
Harpaz has shown that this category is the universal $\infty$-semiadditive category, making the idea mentioned above precise as follows.

\begin{defn}
	We let $\mdef{\Catpifin}$ denote the category of categories admitting $\pi$-finite colimits and functors preserving them.
\end{defn}

\begin{prop}[{\cite[Corollary 5.9]{Harpaz}}]\label{inf-sa-cond}
	A symmetric monoidal category $\CC \in \CAlg(\Catpifin)$ is $\infty$-semiadditive if and only if there is a (necessarily unique) symmetric monoidal $\pi$-finite colimit preserving functor
	\[
		\Span(\Spacespifin) \too \CC.
	\]
\end{prop}

As our proof will rely on the higher descent theorem, which only applies to $\pi$-finite $p$-space colimits, we will need to work with the $p$-typical variant, for which the same proofs show the following.

\begin{defn}
	We let $\mdef{\Catpfin}$ denote the category of categories admitting $\pi$-finite $p$-space colimits and functors preserving them.
\end{defn}

\begin{prop}\label{p-inf-sa-cond}
	A symmetric monoidal category $\CC \in \CAlg(\Catpfin)$ is $p$-typically $\infty$-semiadditive if and only if there is a (necessarily unique) symmetric monoidal $\pi$-finite $p$-space colimit preserving functor
	\[
		\Span(\Spacespfin) \too \CC.
	\]
\end{prop}

Conveniently, it often suffices to check $p$-typical $\infty$-semiadditivity, as the following result shows.

\begin{prop}[{\cite[Corollary 4.4.23]{HL}}]\label{p-loc}
	Let $\CC$ be a $p$-local stable category admitting all limits and colimits.
	If $\CC$ is $p$-typically $\infty$-semiadditive, then it is $\infty$-semiadditive.
\end{prop}

\section{Categorical chromatic localizations}

In this section we recall the chromatic localizations of categories from \cite[\S2]{Desc}, particularly that the $\infty$-semiadditivity of $\SpKn$ implies the $\infty$-semiadditivity of $\CatMn$ (\cref{mn-sa}).

\begin{defn}
	We define the category of \tdef{$\Lnf$-local categories} $\mdef{\CatLnf} \subset \Catperf$ to be the full subcategory on those categories whose mapping spectra are $\Lnf$-local.
	Similarly, we denote $\mdef{\CatLn}$ for the \tdef{$\Ln$-local categories}.
\end{defn}

As explained in \cite[Corollary 2.13]{Desc}, since $\Ln\colon \Sp \to \Ln\Sp$ is smashing localization, we get a smashing localization
\[
	\Ln\colon \Catperf \to \CatLn,
\]
which endows $\CatLn$ with a compatible symmetric monoidal structure.

For any $\CC \in \CatLn$, we can form the fiber
\[
	\Mn(\CC) \too \CC \too \Lnm(\CC),
\]
which is in fact an exact sequence by \cite[Proposition 2.17]{Desc}.
The categories of the form $\Mn(\CC)$ assemble into a category.

\begin{defn}
	We define the category of \mdef{$n$-monochromatic}\footnote{In \emph{op.\ cit.}\ we have used this term for the telescopic variant.} categories $\mdef{\CatMn}$ to be the full subcategory of $\CatLn$ on those $\Ln$-local categories $\CC$ such that $\Lnm(\CC) \simeq \pt$.
\end{defn}

As we have shown in \cite[Proposition 2.17 and Proposition 2.18]{Desc}, the construction $\Mn(\CC)$ appearing above assembles into a coreflection
\[
	\Mn\colon \CatLn \too \CatMn,
\]
which is also smashing in the sense that
\[
	\Mn(\CC) \simeq \CC \otimes \Mn(\Perf(\Ln\SS)),
\]
similarly endowing $\CatMn$ with a compatible symmetric monoidal structure.

Our main interest in $\CatMn$ is its higher semiadditivity.
For the reader's convenience, we provide a sketch of the proof, deviating slightly from the original proof, relying on the equivalent condition from \cref{inf-sa-cond}.

\begin{prop}[{\cite[Theorem 2.29]{Desc}}]\label{mn-sa}
	If $\SpKn$ is $\infty$-semiadditive, then so is $\CatMn$.
\end{prop}

\begin{proof}[Proof sketch]
	Let us denote by $\PrLKn \subset \PrLst$ the full subcategory on those categories whose mapping spectra are $\Kn$-local, and by $\PrLKnw \subset \PrLstw$ their subcategories on the compactly generated categories and functors preserving compact objects (equivalently, whose right adjoint is itself a left adjoint).
	As in \cite[Proposition 2.27]{Desc}, the equivalence
	\[
		\Ind\colon \Catperf \adj \PrLstw \noloc (-)^\omega
	\]
	restricts to an equivalence
	\[
		\Ind\colon \CatMn \adj \PrLKnw \noloc (-)^\omega.
	\]
	Therefore, we equivalently need to show that $\PrLKnw$ is $\infty$-semiadditive.

	Recall that $\PrL$ is $\infty$-semiadditive by \cite[Example 4.3.11]{HL} (see \cite[Proposition 2.21]{CatAmbi} for another proof).
	Thus by \cref{inf-sa-cond} there is a symmetric monoidal $\pi$-finite colimit preserving functor
	\[
		\Span(\Spacespifin) \too \PrL.
	\]
	Recall that $\PrLKn$ is equivalent to $\Mod_{\SpKn}(\PrL)$ by \cite[Proposition 5.2.10]{AmbiHeight}, thus, composing with the symmetric monoidal left adjoint free functor we obtain
	\[
		\Span(\Spacespifin) \too \PrL \too \Mod_{\SpKn}(\PrL) \simeq \PrLKn.
	\]
	Next, we show that this functor factors through the subcategory $\PrLKnw \subset \PrLKn$.
	Observe that the functor in question sends $X$ to $\SpKn[X] \simeq \SpKn^X$, and is given on morphisms by
	\[
		X \oot[f] Z \too[g] Y
		\qquad \mapsto \qquad
		\SpKn^X \too[f^*] \SpKn^Z \too[g_!] \SpKn^Y,
	\]
	see \cite[Theorem B]{CatAmbi} (see also \cite[\S4.1]{CatChar} for a related discussion).
	We note that the right adjoint of $g_!$, namely $g^*$, is itself left adjoint to $g_*$.
	As for $f^*$, its right adjoint is $f_*$.
	Since we assume that $\SpKn$ is $\infty$-semiadditive, $f_*$ is equivalent to $f_!$, and in particular is a left adjoint.
	Finally, recall from \cite[Proposition 2.28]{Desc} that $\PrLKnw \subset \PrLKn$ is a symmetric monoidal subcategory, and is closed under colimits, thus the functor factors into a symmetric monoidal $\pi$-finite colimit preserving functor
	\[
		\Span(\Spacespifin) \too \PrLKnw,
	\]
	which shows that $\PrLKnw$ is $\infty$-semiadditive by \cref{inf-sa-cond}.
\end{proof}

\section{Chromatically localized algebraic K-theory}

In this section we discuss $\Tnp$-localized K-theory.
We start with several recollections, and in \cref{weak-descent} we prove the variant of the higher descent theorem for constant colimits, without higher semiadditivity assumptions.

\begin{defn}
	We let
	\[
		\mdef{\KTnp}\colon \Catperf \too \SpTnp
	\]
	denote the lax symmetric monoidal functor given by the composition $\LTnp\KK(-)$.
\end{defn}

Another key feature of $\CatMn$ is that $\Tnp$-localized K-theory of $\Ln$-local categories factors through monochromatization, due to the vanishing and purity theorems of \cite{purity,DescVan}.

\begin{prop}[{\cite[Proposition 2.24]{Desc}}]\label{purity}
	For any $\Ln$-local category, the map
	\[
		\Mn(\CC) \too \CC
		\qin \CatLn
	\]
	induces an isomorphism
	\[
		\KTnp(\Mn(\CC)) \iso \KTnp(\CC)
		\qin \SpTnp.
	\]
\end{prop}

\begin{proof}
	Consider the exact sequence
	\[
		\Mn(\CC) \too \CC \too \Lnm(\CC).
	\]
	Since K-theory preserves exact sequences, we get an induced exact sequence in $\SpTnp$.
	The last term vanishes by \cite[Theorem C]{DescVan}, so the first map is an isomorphism as required.
\end{proof}

Purity also tells us more about the case of the unit.

\begin{defn}
	We denote by $\mdef{\MM_n} := \Mn(\Perf(\Ln\SS))$ the unit of $\CatMn$.
\end{defn}

\begin{prop}\label{purity-unit}
	There is an isomorphism
	\[
		\KTnp(\MM_n) \iso \KTnp(\SKn)
		\qin \SpTnp.
	\]
\end{prop}

\begin{proof}
	Indeed, by \cite[Purity Theorem]{purity} the map $\Ln\SS \to \SKn$ induces an isomorphism on $\Tnp$-localized K-theory, and the result follows from \cref{purity}.
\end{proof}

We also employ this result to make $\Tnp$-localized K-theory into a lax symmetric monoidal functor out of $\CatMn$.

\begin{prop}\label{lax-catmn}
	The lax symmetric monoidal structure on $\KTnp$ gives a lax symmetric monoidal structure on
	\[
		\KTnp\colon \CatMn \too \SpTnp.
	\]
\end{prop}

\begin{proof}
	Recall that the localization
	\[
		\Ln\colon \Catperf \too \CatLn
	\]
	is compatible with the symmetric monoidal structure, and as a consequence, its right adjoint is lax symmetric monoidal.
	Composing this inclusion with $\KTnp$, we get a lax symmetric monoidal functor
	\[
		\KTnp\colon \CatLn \too \SpTnp.
	\]
	Consider now the colocalization
	\[
		\Mn\colon \CatLn \too \CatMn.
	\]
	This exhibits $\CatMn$ as the Dwyer--Kan localization of $\CatLn$ at the collection of maps $W$ which are sent to isomorphisms under $\Mn$.
	Note that $W$ is closed under tensoring with an object of $\CatLn$, for example since $\Mn$ is smashing.
	Therefore, as explained in \cite[\S3.3.1]{DK}, there is an equivalence between lax symmetric monoidal functors out of $\CatLn[W^{-1}] \simeq \CatMn$ and lax symmetric monoidal functors out of $\CatLn$ sending all morphisms in $W$ to isomorphisms.
	By \cref{purity}, the functor $\KTnp$ indeed sends $W$ to isomorphisms, and the result follows.
\end{proof}

In order to carry out our inductive proof, we shall also make use of the following variant on the higher descent theorem.
We repeat the proof from \emph{op.\ cit.}\ with slight modifications to avoid higher semiadditivity assumptions.

\begin{prop}[{\cite[Theorem 3.1]{Desc}}]\label{weak-descent}
	The functor
	\[
		\KTnp\colon \CatLnf \too \SpTnp
	\]
	preserves constant $\pi$-finite $p$-space colimits.
\end{prop}

\begin{proof}
	We first note that $\KTnp$ commutes with (not necessarily constant) colimits over $1$-finite $p$-spaces.
	Indeed, we may assume that the space is connected, since K-theory commutes with finite products (and the source and the target of $\KTnp$ are $0$-semiadditive).
	For a connected $1$-finite $p$-space, that is, for $A = \BB G$ where $G$ is a finite $p$-group, this is \cite[Theorem C]{DescVan}.
	
	We now need to show that for any $\CC \in \CatLnf$ and any $\pi$-finite $p$-space $A$, the assembly map
	\[
		\KTnp(\CC)[A] \too \KTnp(\CC[A])
	\]
	is an isomorphism.
	
	We reduce to the case where $A$ is simply connected.
	Let $B := \tau_{\leq 1} A$, and let $A_\bullet\colon B \to \Spaces$ denote the straightening of $A \to B$, so that $A = \colim_B A_\bullet$.
	Consider the commutative diagram
	\[\begin{tikzcd}
		{\colim_B \KTnp(\CC)[A_\bullet]} & {\colim_B \KTnp(\CC[A_\bullet])} \\
		& {\KTnp(\colim_B \CC[A_\bullet])} \\
		{\KTnp(\CC)[\colim_B A_\bullet]} & {\KTnp(\CC[\colim_B A_\bullet])}
		\arrow[from=1-1, to=1-2]
		\arrow[from=1-1, to=3-1]
		\arrow[from=1-2, to=2-2]
		\arrow[from=2-2, to=3-2]
		\arrow[from=3-1, to=3-2]
	\end{tikzcd}\]
	Our map is the bottom horizontal map.
	The left vertical map and the lower right vertical map are isomorphisms since the functor $X[-]$ commutes with colimits.
	The space $B$ is $1$-finite, hence the upper right vertical map is an isomorphism by the case already explained above.
	Thus, to show that the bottom horizontal map is an isomorphism, it suffices to show that the top horizontal map is an isomorphism.
	Therefore, it suffices to show that for any $b \in B$, the map
	\[
		\KTnp(\CC)[A_b] \too \KTnp(\CC[A_b])
	\]
	is an isomorphism.
	Since $A_b$ is simply connected, we are reduced to the case where $A$ is simply connected.
	
	We reduce to categories of the form $\CC = \Perf(R)$ where $R$ is an $\Lnf$-local ring spectrum.
	By the Schwede--Shipley theorem \cite{SchwedeShipley} (see also \cite[Proposition 2.9]{Desc}), $\CC$ is a filtered colimit of $\Perf(R)$ for $R$ the endomorphism ring of an object of $\CC$, which is $\Lnf$-local by assumption.
	Thus, by the commutativity of $\KTnp$ with filtered colimits, and the naturality of the assembly map, it suffices to prove the result for $\Perf(R)$, i.e.\ to show that
	\[
		\KTnp(\Perf(R))[A] \too \KTnp(\Perf(R)[A])
	\]
	is an isomorphism.
	Recall that $A$ is (simply) connected.
	Choose an arbitrary base-point, thus by \cite[Corollary A.6]{Desc} the map is identified with
	\[
		\KTnp(R)[A] \too \KTnp(R[\Omega A]).
	\]
	
	We now handle the case where $R$ is of height $0$.
	Since $A$ is simply connected, the map $A \to \pt$ induces an isomorphism
	\[
		\KTo(R)[A] \iso \KTo(R)
	\]
	by \cref{height-1-trunc}.
	Similarly, its loop $\Omega A \to \pt$ induces an isomorphism
	\[
		\KTo(R[\Omega A]) \iso \KTo(R)
	\]
	since $R$ is $p$-invertible while $\Omega A$ is a connected $p$-space, which concludes the case of height $0$.

	We move on to the case where $R$ is of height $n \geq 1$.
	Our proof will proceed by induction on the truncatedness of $A$.
	Write $A$ via the bar construction $A \simeq \colim_{\Delta^\op} A_\bullet$ where $A_k := (\Omega A)^k$.
	Since $A$ is simply connected, hence $\Omega A$ is connected, this induces an isomorphism $\Omega A \simeq \colim_{\Delta^\op} \Omega A_\bullet$.
	Consider the commutative diagram
	\[\begin{tikzcd}
		{\colim_{\Delta^\op}\KTnp(R)[A_\bullet]} & {\colim_{\Delta^\op}\KTnp(R[\Omega A_\bullet])} \\
		& {\KTnp(\colim_{\Delta^\op}R[\Omega A_\bullet])} \\
		{\KTnp(R)[\colim_{\Delta^\op} A_\bullet]} & {\KTnp(R[\Omega\colim_{\Delta^\op} A_\bullet])}
		\arrow[from=1-1, to=1-2]
		\arrow[from=1-1, to=1-2]
		\arrow[from=1-1, to=3-1]
		\arrow[from=1-2, to=2-2]
		\arrow[from=2-2, to=3-2]
		\arrow[from=3-1, to=3-2]
	\end{tikzcd}\]
	Our goal is to show that the bottom horizontal map is an isomorphism, so it suffices to show that all other maps are isomorphisms.
	The top map is an isomorphism by the inductive hypothesis.
	The functor $\KTnp$ preserves sifted colimits by \cite[Corollary 4.32]{purity}, and so do the functors $\Omega(-)$, $R[-]$ and $\KTnp(R)[-]$, so the vertical maps are isomorphisms, concluding the proof.
\end{proof}

\section{The inductive argument}

In this final section we carry out the inductive step of our proof.

\begin{assum}
	The category $\SpKn$ is $\infty$-semiadditive.
\end{assum}

We shall prove in \cref{main-thm} that $\SpKnp$ is $\infty$-semiadditive.

\subsection*{Higher semiadditivity over $\KTnp(\SKn)$}

In this subsection we prove \cref{KTnp-sa}, saying that $\cMod_{\KTnp(\SKn)}$ is $\infty$-semiadditive.

By the inductive assumption and \cref{mn-sa}, the category $\CatMn$ is $\infty$-semiadditive.
Applying \cref{p-inf-sa-cond}, the unit map therefore extends to a map
\[
	\Span(\Spacespfin) \too \CatMn
	\qin \CAlg(\Catpfin),
\]
sending $A$ to $\MM_n[A]$.
Post-compose this with $\KTnp$, which is lax symmetric monoidal by \cref{lax-catmn}.
Since the composition is lax symmetric monoidal, it lifts to modules over the image of the unit, which is $\KTnp(\SKn)$ by \cref{purity-unit}, yielding a lax symmetric monoidal functor
\[
	\Span(\Spacespfin) \too \cMod_{\KTnp(\SKn)}.
\]
In light of \cref{p-inf-sa-cond}, it suffices to show that this functor is strong symmetric monoidal and preserves $\pi$-finite $p$-space colimits.
As an intermediate step, we first prove this for the restriction along the symmetric monoidal functor $\Spacespfin \too \Span(\Spacespfin)$.

\begin{lem}\label{strong-col-spaces}
	The lax symmetric monoidal functor
	\[
		\Spacespfin \too \cMod_{\KTnp(\SKn)}
	\]
	is strong symmetric monoidal and preserves $\pi$-finite $p$-space colimits.
\end{lem}

\begin{proof}
	By \cite[Lemma 2.11]{Harpaz} it preserves $\pi$-finite $p$-space colimits if and only if it preserves \emph{constant} $\pi$-finite $p$-space colimits, and the latter condition is satisfied by \cref{weak-descent}.

	To check that it is symmetric monoidal, first note that it sends the unit to the unit.
	Second, for $A, B \in \Spacespfin$ we need to show that
	\[
		\KTnp(\MM_n[A]) \otimes \KTnp(\MM_n[B]) \too \KTnp(\MM_n[A \times B])
	\]
	is an isomorphism, where $\otimes$ denotes the tensor product in $\cMod_{\KTnp(\SKn)}$.
	Via the preservation of $\pi$-finite $p$-space colimits, this map is equivalent to
	\[
		\KTnp(\MM_n)[A] \otimes \KTnp(\MM_n)[B] \too \KTnp(\MM_n)[A \times B].
	\]
	This map is an isomorphism since the tensor product in $\cMod_{\KTnp(\SKn)}$ commutes with colimits in each variable.
\end{proof}

We now prove the desired properties for the functor itself.

\begin{lem}\label{strong-col-spans}
	The lax symmetric monoidal functor
	\[
		\Span(\Spacespfin) \too \cMod_{\KTnp(\SKn)}
	\]
	is strong symmetric monoidal and preserves $\pi$-finite $p$-space colimits.
\end{lem}

\begin{proof}
	By \cite[Corollary 2.16]{Harpaz}, the functor preserves $\pi$-finite $p$-space colimits if and only if the restriction
	\[
		\Spacespfin \too \Span(\Spacespfin) \too \cMod_{\KTnp(\SKn)}
	\]
	does, which follows from \cref{strong-col-spaces}.
	Similarly, note that the inclusion $\Spacespfin \to \Span(\Spacespfin)$ is symmetric monoidal and essentially surjective, so it suffices to check that the composition is symmetric monoidal, which also follows from \cref{strong-col-spaces}.
\end{proof}

We thus conclude the following.

\begin{thm}\label{KTnp-sa}
	The category $\cMod_{\KTnp(\SKn)}$ is $\infty$-semiadditive.
\end{thm}

\begin{proof}
	Combining \cref{strong-col-spans} with \cref{p-inf-sa-cond}, we get that $\cMod_{\KTnp(\SKn)}$ is $p$-typically $\infty$-semiadditive, and hence $\infty$-semiadditive by \cref{p-loc}.
\end{proof}

\subsection*{Higher semiadditivity over Tate construction}

Our next goal is to pass from modules over $\KTnp(\SKn)$ to modules over $\LKnp(E^{tS^1})$, for any $\Knp$-local commutative ring spectrum $E$.
To that end, we recall the following construction due to Yuan.

\begin{prop}[{\cite[\S3]{yuan}}]\label{yuan-map}
	For any $\Tnp$-local commutative ring spectrum $E$ there is a map of commutative ring spectra
	\[
		\KTnp(E^{tC_p}) \too \LTnp(E^{tS^1})
		\qin \CAlg(\SpTnp).
	\]
\end{prop}

\begin{proof}
	We recall Yuan's construction for the reader's convenience.	
	Since the construction uses trace methods, we first consider the connective cover $e := \tau_{\geq 0} E$.
	Recall that $\THH$ of \emph{commutative} ring spectra is given by colimit over $S^1$.
	More formally, it is the left adjoint of the forgetful $\CAlg(\Sp)^{\BB S^1} \to \CAlg(\Sp)$.
	Thus, for $R \in \CAlg(\Sp)^{\BB S^1}$, the free-forgetful adjunction gives a map
	\[
		\THH(R) \too R
		\qin \CAlg(\Sp)^{\BB S^1},
	\]
	where $\THH(R)$ is the topological Hochschild homology of the underlying commutative ring spectrum of $R$.
	Endowing $e$ with the trivial $S^1$-action, the $C_p$-Tate construction has a residual action by $S^1/C_p \simeq S^1$.
	Applying the above to $R = e^{tC_p}$, we get a map
	\[
		\THH(e^{tC_p}) \too e^{tC_p}
		\qin \CAlg(\Sp)^{\BB S^1}.
	\]
	Taking $S^1$-fixed points, we get
	\[
		\TCm(e^{tC_p}) := \THH(e^{tC_p})^{hS^1} \too (e^{tC_p})^{hS^1/C_p}
		\qin \CAlg(\Sp).
	\]
	Composing with the cyclotomic trace map from $\KK$ to $\TC$ and further to $\TCm$, we obtain
	\[
		\KK(e^{tC_p}) \too (e^{tC_p})^{hS^1/C_p}
		\qin \CAlg(\Sp).
	\]

	By the Tate-orbit lemma \cite[Lemma II.4.2]{NS}, the canonical map
	\[
		e^{tS^1} \too (e^{tC_p})^{hS^1/C_p}
		\qin \CAlg(\Sp)
	\]
	exhibits the target as the $p$-completion of the source, and in particular it induces an isomorphism on $\Tnp$-localizations.
	Composing the inverse of this map with the $\Tnp$-localization of the previous map we get
	\[
		\KTnp(e^{tC_p}) \too \LTnp(e^{tS^1})
		\qin \CAlg(\SpTnp).
	\]

	We now return to $E$ itself.
	By \cite[Corollary 3.2.1]{yuan} (which is a consequence of \cite[Purity Theorem]{purity}), the connective cover map $e \to E$ induces an isomorphism
	\[
		\KTnp(e^{tC_p}) \iso \KTnp(E^{tC_p})
		\qin \CAlg(\SpTnp).
	\]
	We can thus use its inverse to get the desired map
	\[
		\KTnp(E^{tC_p})
		\iso \KTnp(e^{tC_p})
		\too \LTnp(e^{tS^1})
		\too \LTnp(E^{tS^1})
		\qin \CAlg(\SpTnp).
	\]
\end{proof}

We use this to deduce the $\infty$-semiadditivity of $\Tnp$-local modules over $\LKnp(E^{tS^1})$.

\begin{prop}\label{lt-inf}
	The category $\cMod_{\LKnp(E^{tS^1})}$ is $\infty$-semiadditive for any $\Knp$-local commutative ring spectrum $E$.
\end{prop}

\begin{proof}
	Note that it suffices to produce a map of commutative ring spectra from $\KTnp(\SKn)$ to $\LKnp(E^{tS^1})$.
	Indeed, if we have such a map, then base-change along it induces a symmetric monoidal left adjoint functor
	\[
		\cMod_{\KTnp(\SKn)} \too \cMod_{\LKnp(E^{tS^1})}.
	\]
	Since the source is $\infty$-semiadditive by \cref{KTnp-sa}, so is the target by \cref{inf-sa-cond} (or for example by \cite[Corollary 3.3.2(2)]{TeleAmbi}).

	We produce this as a composition of three maps.
	First, the unit $\SKn \to \LKn(E^{tC_p})$ induces
	\[
		\KTnp(\SKn) \too \KTnp(\LKn(E^{tC_p}))
		\qin \CAlg(\SpTnp).
	\]
	Second, by $\Knp$-local Tate vanishing $E^{tC_p}$ is $\Ln$-local, so \cite[Purity Theorem]{purity} implies that $E^{tC_p} \to \LKn(E^{tC_p})$ induces an isomorphism on $\Tnp$-localized K-theory
	\[
		\KTnp(E^{tC_p}) \iso \KTnp(\LKn(E^{tC_p}))
		\qin \CAlg(\SpTnp),
	\]
	and we use its inverse.
	Finally, \cref{yuan-map} gives a map
	\[
		\KTnp(E^{tC_p}) \too \LKnp(E^{tS^1})
		\qin \CAlg(\SpTnp).
	\]
\end{proof}

\subsection*{Higher semiadditivity over $\SKnp$}

Finally, we conclude the proof using a variation on an argument of Hopkins--Lurie \cite[Theorem 5.2.1]{HL} (see in particular \cite[Corollary 5.2.7]{HL}).

\begin{thm}\label{main-thm}
	The category $\SpKnp$ is $\infty$-semiadditive.
\end{thm}

\begin{proof}
	Assume inductively that we have already shown that $\SpKnp$ is $m$-semiadditive, and let $X_\bullet\colon A \to \SpKnp$ be a diagram indexed by an $(m{+}1)$-finite space $A$.
	We wish to show that the norm map associated to this diagram is an isomorphism.
	Observe that the collection of $Y \in \SpKnp$ for which the norm map
	\[
		\Nm\colon \colim (X_\bullet \otimes Y) \too \lim (X_\bullet \otimes Y)
		\qin \SpKnp
	\]
	is an isomorphism forms a thick subcategory of $\SpKnp$.

	We first claim that \cref{lt-inf} shows that the norm map is an isomorphism for $Y = \LKnp(E^{tS^1})$ for any $E \in \CAlg(\SpKnp)$.
	In this case, the diagram $X_\bullet \otimes Y$ lifts to $\cMod_{\LKnp(E^{tS^1})}$, where the associated norm map is an isomorphism by \cref{lt-inf}.
	This implies that the norm map in $\SpKnp$ is an isomorphism as well, since the forgetful functor $\cMod_{\LKnp(E^{tS^1})} \to \SpKnp$ sends the norm map in the source to the norm map in the target (see \cite[Theorem 3.2.3]{TeleAmbi}).
	
	Next, we claim that if $E \in \CAlg(\SpKnp)$ is \emph{complex oriented}, then the norm map is an isomorphism for $Y = E$ as well.
	Note that since $E$ is complex oriented, it is a retract of $E^{tS^1}$.
	Since retracts are preserved under any functor, such as $\Tnp$-localization, $E$ is also a retract of $\LKnp(E^{tS^1})$, and so by thickness the norm is an isomorphism for $Y = E$.

	We apply this to the completed Johnson--Wilson spectrum $\cJWnp := \LKnp\JWnp$.
	We remark that it indeed has a commutative ring spectrum structure by \cite[Theorem 8.3]{cJWEinf} (alternatively, we can apply the claim above to a Lubin--Tate spectrum $\Enp$, of which $\cJWnp$ is a retract).

	Finally, all dualizable $\Knp$-local spectra are in the thick subcategory generated by $\cJWnp$ by \cite[Theorem 8.9]{HSt}.
	In particular, $\SKnp$ is in this thick subcategory, concluding the proof.
\end{proof}

%% file: main.bbl
\begin{thebibliography}{BMCSY24b}

\bibitem[AH68]{AH}
D.W. Anderson and Luke Hodgkin.
\newblock {The K-theory of Eilenberg-Maclane complexes}.
\newblock {\em Topology}, 7(3):317--329, 1968.

\bibitem[BHLS23]{tel}
Robert Burklund, Jeremy Hahn, Ishan Levy, and Tomer~M. Schlank.
\newblock {$K$-theoretic counterexamples to Ravenel's telescope conjecture}.
\newblock 2023.
\newblock \href{https://arxiv.org/abs/2310.17459v1}{arXiv:2310.17459v1}
  [math.AT].

\bibitem[BM24]{CatAmbi}
Shay Ben-Moshe.
\newblock {Categorical Ambidexterity}.
\newblock 2024.
\newblock \href{https://arxiv.org/abs/2411.17281v1}{arXiv:2411.17281v1}
  [math.CT].

\bibitem[BMCSY24a]{card}
Shay Ben-Moshe, Shachar Carmeli, Tomer~M. Schlank, and Lior Yanovski.
\newblock {Chromatic Cardinalities via Redshift}.
\newblock {\em International Mathematics Research Notices},
  2024(14):10918--10924, 05 2024.

\bibitem[BMCSY24b]{Desc}
Shay Ben-Moshe, Shachar Carmeli, Tomer~M. Schlank, and Lior Yanovski.
\newblock {Descent and Cyclotomic Redshift for Chromatically Localized
  Algebraic $K$-theory}.
\newblock {\em Journal of the American Mathematical Society}, 2024.

\bibitem[BR05]{cJWEinf}
Andrew Baker and Birgit Richter.
\newblock {On the $\Gamma$-cohomology of rings of numerical polynomials}.
\newblock {\em Comment. Math. Helv.}, 80(4):691--723, 2005.

\bibitem[BSY24]{null}
Robert Burklund, Tomer~M. Schlank, and Allen Yuan.
\newblock {The Chromatic Nullstellensatz}.
\newblock {\em Annals of Mathematics}, 2024.

\bibitem[CCRY22]{CatChar}
Shachar Carmeli, Bastiaan Cnossen, Maxime Ramzi, and Lior Yanovski.
\newblock {Characters and transfer maps via categorified traces}.
\newblock 2022.
\newblock \href{https://arxiv.org/abs/2210.17364v2}{arXiv:2210.17364v2}
  [math.AT].

\bibitem[CM17]{Tate}
Dustin Clausen and Akhil Mathew.
\newblock {A short proof of telescopic Tate vanishing}.
\newblock {\em Proceedings of the American Mathematical Society},
  145(12):5413--5417, 2017.

\bibitem[CMNN24]{DescVan}
Dustin Clausen, Akhil Mathew, Niko Naumann, and Justin Noel.
\newblock {Descent and vanishing in chromatic algebraic $K$-theory via group
  actions}.
\newblock {\em Annales scientifiques de l'\'Ecole normale sup\'erieure},
  57(4):1135--1190, 2024.

\bibitem[CSY21]{AmbiHeight}
Shachar Carmeli, Tomer~M. Schlank, and Lior Yanovski.
\newblock {Ambidexterity and Height}.
\newblock {\em Advances in Mathematics}, 385:107763, 2021.

\bibitem[CSY22]{TeleAmbi}
Shachar Carmeli, Tomer~M. Schlank, and Lior Yanovski.
\newblock {Ambidexterity in chromatic homotopy theory}.
\newblock {\em Inventiones mathematicae}, 228(3):1145–1254, 2022.

\bibitem[CY23]{CY}
Shachar Carmeli and Allen Yuan.
\newblock {Higher semiadditive Grothendieck-Witt theory and the $K(1)$-local
  sphere}.
\newblock {\em Communications of the American Mathematical Society},
  3(02):65--111, 2023.

\bibitem[GS96]{GS}
John~PC Greenlees and Hal Sadofsky.
\newblock {The Tate spectrum of $v_n$-periodic complex oriented theories}.
\newblock {\em Mathematische Zeitschrift}, 222(3):391--406, 1996.

\bibitem[Har20]{Harpaz}
Yonatan Harpaz.
\newblock {Ambidexterity and the universality of finite spans}.
\newblock {\em Proceedings of the London Mathematical Society},
  121(5):1121--1170, 2020.

\bibitem[Hin16]{DK}
Vladimir Hinich.
\newblock {Dwyer--Kan localization revisited}.
\newblock {\em Homology, Homotopy and Applications}, 18(1):27--48, 2016.

\bibitem[HL13]{HL}
Michael Hopkins and Jacob Lurie.
\newblock {Ambidexterity in $K(n)$-Local Stable Homotopy Theory}.
\newblock \url{https://www.math.ias.edu/~lurie/papers/Ambidexterity.pdf}, 2013.

\bibitem[HS96]{HSa}
Mark Hovey and Hal Sadofsky.
\newblock {Tate cohomology lowers chromatic Bousfield classes}.
\newblock {\em Proceedings of the American Mathematical Society},
  124(11):3579--3585, 1996.

\bibitem[HS99]{HSt}
Mark Hovey and Neil~P Strickland.
\newblock {\em {Morava $K$-theories and localisation}}, volume 666.
\newblock American Mathematical Soc., 1999.

\bibitem[Kuh04]{Kuhn}
Nicholas~J Kuhn.
\newblock {Tate cohomology and periodic localization of polynomial functors}.
\newblock {\em Inventiones mathematicae}, 157(2):345--370, 2004.

\bibitem[LMMT24]{purity}
Markus Land, Akhil Mathew, Lennart Meier, and Georg Tamme.
\newblock {Purity in chromatically localized algebraic $K$-theory}.
\newblock {\em Journal of the American Mathematical Society}, 37:1011--1040,
  2024.

\bibitem[NS18]{NS}
Thomas Nikolaus and Peter Scholze.
\newblock {On topological cyclic homology}.
\newblock {\em Acta Mathematica}, 221(2):203 -- 409, 2018.

\bibitem[RW80]{RW}
Douglas~C. Ravenel and W.~Stephen Wilson.
\newblock {The Morava K-Theories of Eilenberg-MacLane Spaces and the
  Conner-Floyd Conjecture}.
\newblock {\em American Journal of Mathematics}, 102(4):691--748, 1980.

\bibitem[SS03]{SchwedeShipley}
Stefan Schwede and Brooke Shipley.
\newblock {Stable model categories are categories of modules}.
\newblock {\em Topology}, 42(1):103--153, 2003.

\bibitem[Yua24]{yuan}
Allen Yuan.
\newblock {Examples of chromatic redshift in algebraic $K$-theory}.
\newblock {\em Journal of the European Mathematical Society}, 2024.

\end{thebibliography}
